\documentclass[11pt]{article}

\newtheorem{theorem}{Theorem}[section]
\newtheorem{lemma}[theorem]{Lemma}
\newtheorem{prop}[theorem]{Proposition}
\newtheorem{cor}[theorem]{Corollary}
\newtheorem{defn}[theorem]{Definition}
\usepackage{eucal}
\usepackage{amscd}
\usepackage{amssymb}
\usepackage{amsmath}
\usepackage{amsfonts}
\usepackage{amsopn}
\usepackage{latexsym}
\newcommand*{\mx}{\mbox}
\newcommand*{\noin}{{\noindent}}
\newcommand*{\kar}{{\longrightarrow}}
\newcommand*{\mbq}{{\mathbb{Q}}}
\newcommand*{\mt}{{\mathbb{T}}}
\newcommand*{\fp}{{\mathfrak{p}}}
\newcommand*{\fh}{{\mathfrak{h}}}
\newcommand*{\fa}{{\mathfrak{a}}}
\newcommand*{\fq}{{\mathfrak{q}}}
\newcommand*{\fm}{{\mathfrak{m}}}
\newcommand*{\mbr}{{\mathbb{R}}}
\newcommand*{\mbc}{{\mathbb{C}}}
\newcommand*{\mbz}{{\mathbb{Z}}}
\newcommand*{\calok}{{{\cal O}_{K}}}
\begin{document}
\begin{center}
{\bf{\large Ribet's construction of a suitable cusp eigenform}}\\
\vspace{.3cm}
{{\bf Anupam Saikia}}\\
{\it Department of Mathematics, IIT Guwahati,}\\
{\it Guwahati 781039.}\\
{\it Email: a.saikia@iitg.ernet.in}
\end{center}
%-----------------------------------------------------------
\vspace*{.3cm} {\it Abstract}: The aim of this article to give a
self-contained exposition on Ribet's construction of a cusp
eigenform of weight 2 with certain congruence properties for its
eigenvalues.

%----------
%----------------------
\vspace*{.3cm} \noin {\it Acknowledgement}: I am very grateful to
Kevin Buzzard for pointing out errors in an earlier version.
%----------
%----------------------

\section{Preliminaries}
We begin by recalling some of the rudiments of modular forms.
Other basic ingredients are included in the Appendix.
\subsection{Modular forms} Let $p$ be an odd prime. Let
$\mathfrak{h}$ denote the upper half complex plane, i.e.,
\[\mathfrak{h}=\{z\in \mbc \mid Im(z)>0\}.\] Let $SL_{2}(\mbz)$,
$\Gamma_{0}(p)$ and $\Gamma_{1}(p)$ respectively denote the
following groups:
\begin{eqnarray*}
SL_{2}(\mbz)& = & \Bigg\{\left[\begin{array}{cc}
a&b\\c&d\\\end{array}\right]\mid a,\;b,\;c,\;d\;\in
\mbz,\;ad-bc=1\Bigg\}\\
\Gamma_{0}(p)& = & \Bigg\{\left[\begin{array}{cc}
a&b\\c&d\\\end{array}\right]\in SL_{2}(\mbz)\mid \;c\equiv 0
\mbox{ modulo } p \Bigg\}, \\
\Gamma_{1}(p)& = & \Bigg\{\left[\begin{array}{cc}
a&b\\c&d\\\end{array}\right]\in \Gamma_{0}(p)\mid a \equiv 1
\mbox{ modulo } p,\; d\equiv 1\mbox{ modulo } p \Bigg\},
\end{eqnarray*}
Let $GL_{2}(\mbq)$ ($GL_{2}(\mbr)$) denote the $2\times 2$
invertible matrices with rational (real) coefficients. It is easy
to note all these matrix groups act on $\mathfrak{h}$ by sending
$z$ to $\frac{az+b}{cz+d}$.
 For a function $f:
\mathfrak{h}\kar \mbc$ and any fixed integer $k\geq 0$, we can
define a function $f|[\gamma]_{k}$ as
\[f|[\gamma]_{k}(z)=(cz+d)^{-k}f(\gamma(z))\;\;\;\;\forall\;\;\gamma=\left[\begin{array}{cc}
a&b\\c&d\\\end{array}\right]\in GL_{2}(\mbq).\] A function
$f:\mathfrak{h}\kar \mbc$ is called {\it weakly modular} of weight
$k$ with respect to $\Gamma$ if $f|[\gamma]_{k}=f$ for all
$\gamma\in \Gamma$ where $\Gamma$ can mean anyone of
$SL_{2}(\mbz)$, $\Gamma_{0}(p)$ or $\Gamma_{1}(p)$. It is clear
that $\left[\begin{array}{cc} 1&1\\0&1\\\end{array}\right]\in
\Gamma$ and hence we must have $f(z+1)=f(z)$ for a weakly modular
function. If $f$ is holomorphic on $\fh$, we can look at the
Fourier expansion of $f$ in terms of $q=e^{2\pi i z}$, i.e.,
$\sum\limits_{n=-\infty}^{+\infty} a_{n}q^{n}$. We say $f$ is
holomorphic at $\infty$ if its $q$-expansion does not involve
negative powers of $q$, i.e., $a_{n}=0$ for $n<0$. If $a_{n}=0$
for $n\leq 0$, then we say that $f$ vanishes at $\infty$. Note
that $q=e^{2\pi i z}\rightarrow 0$ as $Im(z)\rightarrow \infty$,
justifying the terminology. We say that $f$ is a modular form of
weight $k$ with respect to
$\Gamma$ if \\
(i) $f$ is weakly modular of weight $k$ with respect to $\Gamma$.\\
(ii) $f$ is holomorphic on $\mathfrak{h}$.\\
(iii) $f|[\gamma]_{k}$ is holomorphic at $\infty$ for all
$\gamma\in SL_{2}(\mbz)$. \\
(iv) If, in addition, the $q$-expansion of $f|[\gamma]_{k}$ has
$a(0)=0$ for all $\gamma\in \Gamma$, then $f$ is said to be a
cusp form. \\

\noin Note that it is enough to check the last two conditions for
a finite number of coset representatives $\{\alpha_{i}\}$ of
$\Gamma$ in $SL_{2}(\mbz)$. The set $\{\alpha_{i}(\infty)\}$ is
known as the {\it cusps} of $\Gamma$. Let us denote the space of
all modular forms (cusp forms) of weight $k$ for $\Gamma$ by
$M_{k}(\Gamma)$ ($S_{k}(\Gamma)$ respectively). These turn out to
be finite dimensional vector spaces. The quotient vector space of
$M_{k}(\Gamma)$ by $S_{k}(\Gamma)$ is known as the Eisenstein
space, denoted by ${\cal E}_{k}(\Gamma)$. It can be identified as
the orthogonal complement of $S_{k}(\Gamma)$ under Petersson inner
product, and hence can be thought of as a subspace of
$M_{k}(\Gamma)$ (see section 6.6 of Appendix).

\subsection{Semi-cusp forms}

\begin{defn}A semi-cusp form $f$ is a modular form whose leading Fourier
coefficient is $0$, though $f|[\gamma]_{k}$ need not have its
leading Fourier coefficient $0$ for all $\gamma\in SL_{2}(\mbz)$.
In other words, a semi-cusp form vanishes at $\infty$, but it need
not vanish at the other `cusps'. We shall denote the space of
semi-cusp forms of $\Gamma$ by $S'_{k}(\Gamma)$.
\end{defn}

\noin Consider the map
\[\beta: \Gamma_{0}(p)\kar (\mbz/p\mbz)^{\times}, \;\;\gamma=\left[\begin{array}{cc}
a&b\\c&d\\\end{array}\right]\mapsto d\mx{ mod }p.\]  (Note that
$(d,p)=1$ for $\gamma\in \Gamma_{0}(p)$ as $ad-bc=1$ and $p|c$).
Clearly, $\Gamma_{1}(p)$ is the kernel of $\beta$, and the
quotient is $(\mbz/p\mbz)^{\times}$.
%Thus, $\Gamma_{0}(p)$ acts on $M_{k}(\Gamma_{1}(p))$ through
%the quotient $(\mbz/p\mbz)^{\times}$.
For a character $\epsilon$ of $(\frac{\mbz}{p\mbz})^{\times}$, we
can define a subspace $M_{k}(\Gamma_{1}(p), \epsilon)$ of
$M_{k}(\Gamma_{1}(p))$, which consists of modular forms $f$ such
that $f|[\gamma]_{k}=\epsilon(d)f$ for any
$\gamma=\left[\begin{array}{cc}
a&b\\c&d\\\end{array}\right]\in\Gamma_{0}(p)$. We can define
$S'_{k}(\Gamma_{1}(p), \epsilon)$ and $S_{k}(\Gamma_{1}(p),
\epsilon)$ analogously. Note that any character of
$\big(\frac{\mbz}{p\mbz}\big)^{\times}$ is of the form $w^{i}$,
$i=0,\;1,\;\ldots,\; (p-2)$ where $w$ is the Teichmuller character
(see section 6.5 Appendix).

%----------
%----------------------

\subsection{Examples of modular forms}

For a non-trivial even character $\epsilon$ of
$(\frac{\mbz}{p\mbz})^{\times}$,  we have the following Eisenstein
series of weight 2 and type $\epsilon$ (cf chapter 4 of [Di-S]:
\begin{eqnarray}
\label{e1}G_{2,\epsilon}& = & \frac{L(-1,\epsilon)}{2}+
\sum_{n\geq
1}\sum_{d|n}\epsilon(d)d q^{n},\\
\label{e2} s_{2,\epsilon}& = &\sum_{n\geq
1}\sum_{d|n}\epsilon\big(\frac{n}{d}\big)d q^{n}.\end{eqnarray}

\noin These two form a basis for the Eisenstein space ${\cal
E}_{2}(\Gamma_{1}(p),\epsilon)$ (cf theorem 4.6.2 [Di-S]). Note
that $s_{2,\epsilon}$ is a semi-cusp form. Moreover, both of these
are eigenvectors for all Hecke operators $T_{l}$ with $(l,p)=1$
(cf proposition 5.2.3 [Di-S]):
\begin{equation*}
T_{l}s_{2,\epsilon}=(l+\epsilon(l))s_{2,\epsilon}, \qquad
T_{l}G_{2,\epsilon}=(1+\epsilon(l)l)G_{2,\epsilon}.\] (See section
6.7 of the Appendix for Hecke operators.)
% These two forms have coefficients defined over the ring $\mbz[\mu_{p-1}]$ of
%integers in $K=\mbq(\mu_{p-1})$.

If $\epsilon$ is an odd character of
$\big(\frac{\mbz}{p\mbz}\big)^{\times}$, we have an Eisenstein
series of weight 1 and type $\epsilon$ given by (cf section 4.8 in
[Di-S])
\[G_{1,\epsilon}=\frac{L(0,\epsilon)}{2}+ \sum_{n\geq 1}\sum_{d|n}\epsilon(d)
q^{n}.\] The above three forms have coefficients defined over
$\mbq(\mu_{p-1})$, where $\mu_{p-1}$ denotes the $(p-1)^{th}$
roots of 1. Let $\wp$ denote any of the unramified primes of
$\mbq(\mu_{p-1})$ lying above $p$. Clearly, all the Eisenstein
forms given above have $\wp$ integral coefficients (except
possibly for the constant terms, but see lemma 3.1 later).

For the trivial character $\epsilon=1$, we have the following
Eisenstein series (cf Theorem 4.6.2 in [Di-S]) in
$M_{k}(\Gamma_{0}(p))=M_{k}(\Gamma_{1}(p),1)$:
\begin{eqnarray}
\label{eis1} G_{k}& = & -\frac{B_{k}}{2k}+ \sum_{n\geq
1}\sum_{d|n}d^{k-1} q^{n} \mx{ for }k\geq 4,\\
G_{2}& = & E_{2}(z)-pE_{2}(pz),\mx{ where
}E_{2}(z)=-\frac{B_{2}}{4}+ \sum_{n\geq 1}\sum_{d|n}dq^{n},
\end{eqnarray}

%----------
%----------------------

\section{Key steps in the construction of the unramified
$p$-extension}

For Ribet's construction of an unramified extension of
$\mbq(\mu_{p})$, one requires a Galois representation on which the
Frobenius elements act in a suitable way (see[D]). We can use the
representation associated with a cusp eigenform (cf chapter 9 of
[Di-S]). But we need to show that there indeed exists a cusp
eigenform whose eigenvalues have certain congruence properties.

The Eisenstein series $G_{2,\epsilon}$ is a simultaneous eigenform
for the Hecke operators $T_{l}$ where $l$ is a prime other than
$p$, with corresponding eigenvalues $1+\epsilon(l)l\equiv
1+l^{k-1}$ modulo $\wp$. Here, $\wp$ denotes a prime of
$\mbq(\mu_{p-1})$ lying above $p$. It turns out that we need
precisely these congruence properties for the Hecke eigenvalues of
a {\it cusp} form. Ribet's idea is to subtract off the constant
term of the Eisenstein series $G_{2,\epsilon}$ in a way that
preserves the congruence properties of the coefficients and leaves
us with a semi-cusp form $f$ which is an eigenvector modulo $\wp$
for all Hecke operators $T_{l}$ with $(l,p)=1$. Then one can
invoke a result of Deligne and Serre and obtain a semi-cusp form
$f'$ which is also an eigenvector for the $T_{l}$'s with
eigenvalues congruent to those of $f$ modulo $\wp$. The congruence
properties of $f'$ then ensures that $f'$ is actually a cusp form.
Any cusp form in $S_{2}(\Gamma_{1}(p))$ is bound to be a newform.
Thus, one can invoke the theory of newforms to conclude that $f'$
is in fact a cusp eigenform, that is, an eigenvector for all Hecke
operators including $T_{n}$'s with $p|n$.

To remove the constant term of the Eisenstein series
$G_{2,\epsilon}$ without affecting the congruence properties of
its coefficients modulo $\wp$, it suffices to produce another
Eisenstein series whose constant term is a $\wp$-unit. This will
be done in the next section.

\section{Construction of an Eisenstein series with $\wp$-unit constant term}

As before, we will denote by $\wp$ a prime of $\mbq(\mu_{p-1})$
lying above $p$. Note that $\wp$ is unramified. We continue to
denote the Teichmuller character by $w$.

\begin{lemma} Let $k$ be even and $2\leq k \leq p-3$. Then the
$q$-expansions of the modular forms $G_{2,w^{k-2}}$ and
$G_{1,w^{k-1}}$ have $\wp$-integral coefficients in
$\mbq(\mu_{p-1})$ and are congruent modulo $\wp$ to the
$q$-expansion
\[-\frac{B_{k}}{2k}+ \sum_{n\geq 1}\sum_{d|n}d^{k-1}
q^{n}.\]

\end{lemma}
Proof: Since $w(d)\equiv d$ mod $\wp$, $w^{k-2}(d)d\equiv d^{k-1}$
mod $\wp$ and $w^{k-1}(d)\equiv d^{k-1}$ mod $\fp$. Hence it
suffices to investigate the constant terms only. We know that (see
(\ref{l1}) and (\ref{l2}) of Appendix)
\begin{eqnarray*}
L(0,\epsilon)& = &
\frac{-1}{p}\sum_{n=1}^{p-1}\epsilon(n)\Big(n-\frac{p}{2}\Big),\\
L(-1,\epsilon)& = &
\frac{-1}{2p}\sum_{n=1}^{p-1}\epsilon(n)\Big(n^{2}-pn-\frac{p^{2}}{6}\Big).
\end{eqnarray*}
Since we know that $w(n)\equiv n^{p}$ mod ($\wp^{2}$) (cf section
6.5 of Appendix), we find that \begin{eqnarray*} pL(0,w^{k-1})&
\equiv &
-\sum_{n=1}^{p-1}n^{1+p(k-1)}\mx{ mod }\wp^{2},\\
pL(-1,w^{k-2})& \equiv & -\frac{1}{2}\sum_{n=1}^{p-1}n^{2+p(k-2)}
\mx{ mod }\wp^{2}.\end{eqnarray*} Note that
$\sum\limits_{n=1}^{p-1}\epsilon(n)n\equiv 0$ mod $\wp$ when
$\epsilon$ is an even character. Moreover, we know that (see
proposition 6.6 of Appendix)
\[pB_{t}\equiv \sum_{n=1}^{p-1}n^{t}\mx{ mod
}p^{2}.\] Therefore, we have
\begin{eqnarray*}
 L(0,w^{k-1})&
\equiv & -\frac{1}{2}B_{1+p(k-1)} \equiv
-\frac{1}{2}(1+p(k-1))\frac{B_{k}}{k}\equiv -\frac{B_{k}}{k} \mx{
mod }\wp, \\ L(-1,w^{k-2})& \equiv & -\frac{1}{2}B_{2+p(k-2)}
\equiv
-\frac{1}{2}(2+p(k-2))\frac{B_{k}}{k}\equiv  -\frac{B_{k}}{k} \mx{ mod }\wp.\\
\end{eqnarray*}
For the second equivalence of each statement above, we use Kummer
congruence as explained in proposition 6.4 in the Appendix.
 Note that
\begin{eqnarray*}
1+p(k-1)& = & k+(p-1)(k-1)\equiv k \mx{ mod }(p-1),\\
2+p(k-2)& = & k+(p-1)(k-2)\equiv k\mx{ mod }(p-1).  \qquad\square
\end{eqnarray*}

\noin The following corollary is now obvious.
\begin{cor}
Let $k$ be even and $2\leq k \leq p-3$. Let $n$, $m$ be even
integers such that $n+m \equiv k$ mod $(p-1)$ and $2\leq n,\;m\leq
p-3$. The the product $G_{1,w^{n-1}}G_{1,w^{m-1}}$ is a modular
form of weight 2 and type $w^{k-2}$ whose $q$-expansion
coefficients are $\wp$-integral in $\mbq(\mu_{p-1})$. Its constant
term is a $\wp$-adic unit if neither $B_{n}$ nor $B_{m}$ is
divisible by $p$.
\end{cor}

\noin The next theorem guarantees the existence of the Eisenstein
series we are looking for.

\begin{theorem}
Let $k$ be an even integer $2\leq k\leq p-3$. Then there exists a
modular form $g$ of weight 2 and type $w^{k-2}$ whose
$q$-expansion coefficients are $\wp$-integers in $\mbq(\mu_{p-1})$
and whose constant term is a $\wp$-unit.
\end{theorem}
Proof:\\
{\underline{Case (i)}} If $p\not|B_{k}$, we can take $G_{2,w^{k-2}}$ by
lemma 3.1. \\
{\underline{Case (ii)}} If we have a pair of even integers $m\;n$
such that $n+m \equiv k$ mod $(p-1)$, $2\leq n,\;m\leq p-3$ and
$p\not|B_{m}B_{n}$, then we can take $G_{1,w^{n-1}}G_{1,w^{m-1}}$
by corollary 3.2.\\
{\underline{Case (iii)}} Suppose neither of the above two cases
are true. We will show that consequently too many Bernoulli
numbers will be $p$-divisible, which will lead to violation of an
upper bound for the $p$-part $h_{p}^{*}$ of the relative class
number of $\mbq(\mu_{p})$. Let $t$ be the number of even integers
$n$, $2\leq n\leq p-3$ such that $p$ divides $B_{n}$. It is easy
to see that $t\geq \frac{p-1}{4}$ if the cases (i) and (ii) do not
arise. But then, $p^{t}$ must divide $h_{p}^{*}$ (see section 6.2
of Appendix). However, that contradicts a result of Carlitz, which
says that $h_{p}^{*}<p^{(\frac{p-1}{4})}$. Hence we must be in
either in case (i) or case (ii). \hspace*{1cm}$\square$

%----------
%----------------------

\section{Existence of a semi-cusp form with suitable eigenvalues}

In this section, we will first construct a semi-cusp form $f$
which is a simultaneous eigenvector modulo $\wp$ for all Hecke
operators $T_{l}$ with $(p,l)=1$. Then we will lift $f$ to a
semi-cusp form $f'$ which is an eigenvector for all such
$T_{l}$'s.

Fix an even integer $k$, $2\leq k\leq p-3$ and assume that
$p|B_{k}$. Consider $\epsilon =w^{k-2}$. Since
$B_{2}=\frac{1}{6}$, $k$ is at least $4$, and hence $\epsilon$ is
a non-trivial even character. We will only be interested in
modular forms of weight 2 and type $\epsilon$.

%----------
%----------------------
\begin{prop}
There exists a semi-cusp form $f=\sum\limits_{n\geq 1}a_{n}q^{n}$
such that $a_{n}$ are $\wp$-integers in $\mbq(\mu_{p-1})$ and such
that $f\equiv G_{2,\epsilon}\equiv G_{k}$ mod $\wp$.
\end{prop}
Proof: Consider $f=G_{2,\epsilon}-c.g$, where $c$ is the constant
term of $G_{2,\epsilon}$. Then $f$ is a semi-cusp form. Now, $c\in
\wp$ as $p|B_{k}$. Hence, $f\equiv G_{2,\epsilon}\equiv
G_{k}$ mod $\wp$. \hspace*{1cm}$\square$\\

\noin Observe further that $f$ is a mod $\wp$-eigenform for all
Hecke operators $T_{l}$ with $(l,p)=1$, as the Eisenstein series
$G_{2,\epsilon}$ is an eigenform form for all such $T_{l}$ with
eigenvalue $(1+\epsilon(l)l)$. Therefore,
\begin{equation}
\label{modp} T_{l}(f)\equiv T_{l}(G_{2,\epsilon})\equiv
(1+\epsilon(l)l)G_{2,\epsilon}\equiv (1+\epsilon(l)l)f \mx{ modulo
} \wp.
\end{equation}

%----------
%----------------------

\subsection{Deligne-Serre lifting lemma}

The following result of Deligne and Serre [D-S] ensures that there
exists a semi-cusp form $f'$ which is an eigenvector for the
$T_{l}$'s ($(l,p)=1$) with eigenvalues congruent modulo $\wp$ to
those of the mod-$\wp$ eigenvector $f$ obtained previously.
\begin{lemma}
Let $M$ be a free module of finite rank over a discrete valuation
ring $R$ with residue field $k$, fraction field $K$ and maximal
ideal $\mathfrak{m}$. Let $S$ be a (possibly infinite) set of
commuting $R$-endomorphisms of $M$. Let $0\neq f\in M$ be an
eigenvector modulo $\mathfrak{m}M$ for all operators in $S$, i.e.,
$Tf=a_{T}f$ mod $\mathfrak{m}M$ $\forall T\in S$ ($a_{T}\in R$).
Then there exists a DVR $R'$ containing $R$ with maximal ideal
$\mathfrak{m}'$ containing $\mathfrak{m}$, whose field of
fractions $K'$ is a finite extension of $K$ and a non-zero vector
$f'\in R'\otimes_{R}M$ such that $Tf'=a'_{T}f'$ for all $T\in S$
with eigenvalues $a'_{T}$ satisfying $a'_{T}\equiv a_{T}$ mod
$\mathfrak{m}'$.
\end{lemma}
Proof: Let $\mt$ be the algebra generated by $S$ over $R$. Clearly
$\mt\in End_{R}(M)$. As $M$ is an free $R$-module of finite rank,
so is $End_{R}(M)$. Therefore, $\mt$ is also free module of finite
rank over $R$, generated by $T_{1},\;\ldots ,\;T_{r}\in S$. Let
$h_{i}$ denote the minimal polynomial of $T_{i}$ acting on
$K\otimes_{R} M$. If we adjoin the roots of all such minimal
polynomials to $K$, we get a finite extension $K'$ of $K$. The
integral closure of $R$ in $K'$ gives us a DVR $R'$ with maximal
ideal ${\fm}'$ lying over $m$, and with residue field $k'$
containing $k$. By replacing $M$ with $R'\otimes M$ and $\mt$ with
$R'\otimes_{R} \mt$, we will continue to write $R$, $\fm$, $k$,
$K$ in stead of $R'$, $\fm$ etc.

Consider the ring homomorphism $\lambda:\mt \kar k$ given by
$T\mapsto a_{T}$ mod $\fm$ for all $T$ in $S$. Clearly,
$ker(\lambda)$ is a maximal ideal of $\mt$. Choose a minimal prime
$\wp$ in $ker(\lambda)$. Then, $\wp$ is contained in the set of
zero-divisors of $\mt$ (see proposition 6.9 of Appendix). As $\mt$
is a free $R$-module, $R$ contains no zero-divisors of $\mt$ and
hence, $\mathfrak{p}\cap R=\{0\}$. Thus, $\mt/\mathfrak{p}$ is a
finite integral extension of $R$. Let $L$ denote the field of
fractions of the integral domain $\mt/\mathfrak{p}$. Let $R_{L}$
be the integral closure of $R$ in $L$, then $R_{L}$ is a DVR with
maximal ideal $m_{L}$ containing $\fm$ and residue field $l$
containing $k$.

Consider the map $\lambda':\mt \kar
\mt/\mathfrak{p}(\hookrightarrow R_{L})$ given by reduction modulo
$\mathfrak{p}$. Let $\lambda'(T)=a'_{T}$ for all $T\in S$.
Clearly, $\lambda'$ maps the maximal ideal $ker(\lambda)$ of $\mt$
into the maximal ideal $m_{L}$ of $R_{L}$. But $(T-a_{T})\in ker
(\lambda)$, hence $\lambda'(T-a_{T})\in m_{L}$ i.e., $a'_{T}\equiv
a_{T}$ modulo $m_{L}$.

Now consider the ring $K\otimes_{R} \mt$. It is an Artinian ring,
hence it has finitely many maximal ideals with residue fields all
isomorphic to $K$. Let ${\cal P}$ be the prime ideal in
$K\otimes\mt$ generated by $\mathfrak{p}$. It will suffice to show
that ${\cal P}$ is an associated prime of $K\otimes M$. Note that
$\wp\subset ker(\lambda)$ implies $\wp$ annihilates $f$ in
$M/\fm$. Now let $x\in Ann_{\mt/\fm}(f)$, say
$x=\bar{g}(T_{1},\ldots, T_{n})$. Then, $x=\bar{g}(a'_{T_{1}},
\ldots,a'_{T_{1}})$ modulo $(T_{1}-a'_{T_{1}}, \ldots,
T_{n}-a'_{T_{n}})$. Thus, $xf=\bar{g}(a'_{T_{1}},
\ldots,a'_{T_{1}})f$ modulo $m_{L} M$, noting that $T-a'_{T}\in
\wp$, and $\wp$ annihilates $f$ modulo $m_{L} M$. As $a'_{T}\equiv
a_{T}$ mod $m_{L}$, we must have $\bar{g}(a_{T_{1}},
\ldots,a_{T_{1}})f=0$ mod $m_{L} M$. As $f\neq 0$, we must have
$\bar{g}(a_{T_{1}}, \ldots,a_{T_{1}})=0$ in $l$. Thus, $x\in \wp$,
and $\wp=Ann_{\mt/\fm}(f)$ is an associated prime of $M/\fm$. For
proof of the following two statements,
see section 6.8.2 of Appendix. \\
(i) $\mathfrak{p}$ is in $Assoc_{\mt/m}(M/\fm)$, hence in
$Supp_{\mt/m}(M/\fm)$, and hence $Ann_{\mt/m}(M/\fm)\subset \wp$. \\
(ii) Now, it follows that $Ann_{K\otimes \mt}(K\otimes M)\subset
{\cal P}$, hence ${\cal P}\in Supp_{K\otimes \mt}(K\otimes M)$ and
therefore ${\cal
P}$ is in $Assoc_{K\otimes \mt}(K \otimes M)$.\\

\noin Now, ${\cal P}$ is the annihilator of some $0\neq f''\in
K\otimes M$, hence ${\cal P}$ annihilates some $f'\in M$. As
$T-a'_{T}\in \mathfrak{p}$, we have $T-a'_{T}\in {\cal P}$ and
$(T-a'_{T})(f')=0$. Thus, $Tf'=a'_{T}f'$ where $a'_{T}\equiv
a_{T}$ modulo $m_{L}$, which concludes our proof.
\hspace*{1cm}$\square$\\

\subsection{Lifting the semi-cusp form to an eigenvector for
$T_{n}$ for $(n,p)=1$}

The following theorem ensures that we have a semi-cusp form which
is an eigenvector for all Hecke operators $T_{n}$ with $p\not|n$.

\begin{theorem}
There is a semi-cusp form $f'=\sum_{n=1}^{\infty}c_{n}q^{n}$ of
weight 2 and type $\epsilon$ such that all its coefficients are
defined over a finite extension of $L$ of $\mbq(\mu_{p-1})$ and
are $\wp_{L}$-integral where $\wp_{L}$ is a prime above $p$.
Further, $T_{l}f'\equiv (1+\epsilon(l)l)f'$ modulo $\wp_{L}$.

\end{theorem}
Proof: There is a basis $B$ of $S'_{2}(\Gamma_{1}(p),\epsilon)$
consisting of semi-cusp forms all of whose coefficients are
defined over a finite extension $K$ of $\mbq(\mu_{p-1})$. Let $R$
be the localization of the ring of integers of $K$ at a prime
$\wp_{K}$ above $\wp$. Let $M$ be the free $R$-module of semi-cusp
forms generated by $B$. Let $S=\{T_{n}|(p,n)=1\}$. We know by
proposition 4.1 and (\ref{modp}) that there exists $f\in M$ such
that
\[T_{l}(f)\equiv (1+\epsilon(l)l)f \mx{ modulo } \wp.\]
By applying the lifting lemma 4.2, we can conclude that there is a
finite extension $L$ of $K$ with a prime $\wp_{L}$ over $\wp_{K}$
such that there exists a semi-cusp form $f'$, with
$\wp_{L}$-integral coefficients in $L$ such that
$T_{l}(f')=c_{l}f'$ and $c_{l}\equiv 1+\epsilon(l)l$ modulo
$\wp_{L}$. \hspace*{1cm} $\square$

%----------
%----------------------

\section{Construction of cusp eigenform}

We will first show that the semi-cusp form $f'$ obtained in the
previous section is in fact a cusp form. Then, we will finally
show that the cusp form $f'$ must be an eigenvector for all Hecke
operators $T_{n}$ including those $n$ which are not co-prime to
$p$.

\subsection{Existence of a suitable cusp form}

\begin{prop}
There exists a non-zero cusp form $f'$ of type $\epsilon$, which
is an eigenform for all Hecke operators $T_{n}$ with $(n,p)=1$,
and which has the property that for any prime $l\neq p$, the
eigenvalue $\lambda_{l}$ of $T_{l}$ acting on $f'$ satisfies
\[\lambda_{l}\equiv1+l^{k-1}\equiv 1+\epsilon(l)l\mx{ mod } {\wp_{L}},\]
where $\wp_{L}$ is a certain prime (independent of $l$) lying
over $\wp$ in the field $L=\mbq(\mu_{p-1}, \lambda_{n})$ generated
by the eigenvalues over $\mbq(\mu_{p-1})$.
\end{prop}
Proof: We already established the existence of a semi-cusp form
$f'$ which is an eigenform for all Hecke operators $T_{n}$
$(n,p)=1$ whose eigenvalues have the required congruence
properties. It suffices to assert that $f'$ is in fact a cusp
form. As $M_{2}(\Gamma_{0}(p), \epsilon)$ is spanned by the cusp
forms, the semi-cusp form $S_{2,\epsilon}$ and the Eisenstein
series $G_{2,\epsilon}$, we must have
\[S'_{2}(\Gamma_{1}(p), \epsilon)=S_{2}(\Gamma_{1}(p),
\epsilon)\oplus \mbc s_{2,\epsilon},\] where orthogonality of the
Eisenstein space and the space of cusp forms under Petersson inner
product $<,>$ is the reason behind the above sum being a direct
one (see section 6.6 of Appendix). Suppose $f'=h+as_{2,\epsilon}$
($a\neq 0$). Then, $f'-as_{2,\epsilon}\in S_{2}(\Gamma_{1}(p),
\epsilon)$. But, $f'-as_{2,\epsilon}\in {\cal
E}_{2}(\Gamma_{1}(p), \epsilon)$ as well, where ${\cal
E}_{2}(\Gamma_{1}(p), \epsilon)$ denotes the subspace consisting
of Eisenstein series in $M_{2}(\Gamma_{1}(p),\epsilon)$. As the
orthogonal subspaces ${\cal E}_{2}(\Gamma_{1}(p), \epsilon)$ and
$S_{2}(\Gamma_{1}(p), \epsilon)$ have trivial intersection,
$f'-as_{2,\epsilon}=0$, i.e., $f'=as_{2,\epsilon}$. Applying
$T_{l}$ to both sides, $(l\neq p)$, we see that we must have
$1+\epsilon(l)l\equiv l+\epsilon(l)$ mod $\wp_{L}$, which forces
$\epsilon(l)=1$. But $\epsilon$ is a non-trivial character and
$l\neq p$ is arbitrary, hence $f'$ must be a cusp form. $\square$

%----------
%----------------------

\subsection{Operators $T_{n}$ for $(n,p)\neq 1$}

So far, we know that we have a cusp form $f$ for $\Gamma_{1}(p)$
of weight 2 and type $\epsilon$ which is an eigenform for all
Hecke operators $T_{l}$ $(l,p)=1$. In this section we will assert
that $f$ is in fact a common eigenform for all Hecke operators,
including $T_{n}$ $(n,p)\neq 1$.

\begin{prop}
Any form $f'$ as above is an eigenform for all Hecke operators
(including those for which $p|n$). Hence, after replacing $f'$ by
a suitable multiple of $f'$, we have
\[f'=\sum_{n=1}^{\infty}\lambda_{n}q^{n}, \mx{ where }
T_{n}(f')=\lambda_{n}f'.\]\end{prop}

Proof: $f'$ must be a newform. For, if it were an old form it will
have to originate from a non-zero modular form in
$M_{2}(SL_{2}(\mbz))$, but that space is trivial. Now for a new
form $f'$, if it is an eigenform for $T_{n}$ ($(n,p)=1$) it has to
be an eigenform for all $T_{n}$ by the theory of newforms (see
Theorem 5.8.2 of [Di-S]). Now we can take a suitable multiple of
$f'$ to get a normalized cusp eigenform as prescribed in the
theorem. \hspace*{1cm}$\square$\\

\noin {\it{Remark}}: The cusp eigenform obtained above can be
associated to a Galois representation which finally gives an
unramified $p$-extension of $\mathbb{Q}(\mu_{p})$, where $\mu_{p}$
denotes the $p$-power roots of unity for an odd prime $p$. This
exposition can be found in [D].

%----------
%----------------------

\section{Appendix}
Here we provide a brief discussion of the various ingredients used
in the previous sections.
%----------
%----------------------

\subsection{Dirichlet $L$-functions}

A Dirichlet character is a homomorphism $\chi:
\big(\frac{\mbz}{N\mbz}\big)^{\times}\kar \mbc^{\times}$, where
$N$ is any positive integer, and $A^\times$ denote the
multiplicative group of units in a ring $A$. $N$ is called the
conductor of $\chi$ if $\chi$ does not factor through
$\big(\frac{\mbz}{M\mbz}\big)^{\times}$ for any $M<N$. We denote
the conductor of $\chi$ by $f_{\chi}$. We can easily extend the
definition of $\chi$ to $\mbz$ by setting $\chi(n)=\chi(n \mx{ mod
}N)$ if $(n,N)=1$ and $\chi(n)=0$ otherwise. The Dirichlet
$L$-series of $\chi$ is defined as
\[L(s,\chi)=\sum_{n=1}^{\infty}\chi(n)n^{-s},\]
where $s$ is a complex number with $Re(s)>1$. It is well-known
that $L(s,\chi)$ can be analytically continued to the whole
complex plane except a simple pole of residue 1 at $s=1$ when
$\chi$ is the trivial character (in which case the function is
just the Riemann-zeta function). Further, $L(s,\chi)$ satisfies a
functional equation relating its values at $s=1$ to values $1-s$.
It also has a Euler product, i.e.,
\[L(s,\chi)=\prod_{l}(1-\chi(l)l^{-s})^{-1}, \;\;Re(s)>1\]
where $l$ runs over the rational primes. The Dirichlet
$L$-functions are related to the Dedekind zeta function of an
abelian number field, as explained below.

Recall that for a number field $K$, the Dedekind zeta function is
defined as
\[\zeta_{K}(s)=\sum_{\fa}(N\fa)^{-s},\;\;\;\;Re(s)>1,\]
where $\fa$ runs over the ideals of the ring $\calok$ of integers
in $K$. It is well-known that $\zeta_{K}(s)$ can be analytically
continued to the whole complex plane except for a simple pole at
$s=1$. Further, $\zeta_{K}(s)$ satisfies a functional equation,
relating the values at $s$ to values at $1-s$.

We can view $\chi$ as a Galois character
\[ \chi \;:\; Gal(\mbq(\mu_{N})/\mbq)\simeq (\mbz/N\mbz)^{\times}\kar
\mbc^{\times},\] and this gives a correspondence $\chi\rightarrow
$ fixed subfield of $ker(\chi)$ in $\mbq(\mu_{N})$, which is an
abelian extension of $\mbq$. This leads to a one-to-one
correspondence between groups of Dirichlet characters and abelian
extensions of $\mbq$. If $K$ is an abelian extension of $\mbq$, it
is contained in some $\mbq(\mu_{N})$ and there will be a
corresponding group $X$ of Dirichlet characters of conductor
dividing $N$.

If $K$ is an abelian number field and $X$ is the corresponding
group of Dirichlet characters, then one can show that (see theorem
4.3 in [Wa])
\[\zeta_{K}(s)=\prod_{\chi \in X}L(s,\chi).\]

%----------
%----------------------

\subsection{The relative class number and Dirichlet $L$-values}

 The analytic class number formula is given by
\[\lim_{s\rightarrow
1}\zeta_{K}(s)=\frac{2^{r_{K}}(2\pi)^{t_{K}}h_{K}R_{K}}{w_{K}\sqrt{|d_{K}|}},\]
where $r_{K}$ and $t_{K}$ denote respectively the number of real
and complex pairs of embedding of $K$, $w_{K}$ the number of roots
of unity in $K$, $R_{K}$ the regulator of $K$, $d_{K}$ the
discriminant of $K$ and $h_{K}$ the class number of $K$.

Now consider $K=\mbq(\zeta_{p})$, then $r_{K}=0$,
$t_{K}=\frac{p-1}{2}$. Let $K^{+}$ be the maximal real subfield of
$K$, for which $r_{K^{+}}=\frac{p-1}{2}$ and $t_{K^{+}}=0$. It is
easy to establish that $h_{K^{+}}$ divides $h_{K}$. The relative
class number of $K$ is defined as
$h_{K}^{-}=\frac{h_{K}}{h_{K^{+}}}$. The purpose of this section
is to investigate the $p$-part $h_{K}^{-}$, and relate it to the
values of Dirichlet $L$-functions.

\begin{prop}
\[h_{K}^{-}=\alpha p\prod_{i=0}^{p-2}L(0,w^{i}),\]
where $\alpha$ is a certain power of $2$.
\end{prop}
Proof: Dividing the analytic class number formulas for $K$ and
$K^{+}$, and then shifting the limit to $s\rightarrow 0$ via the
functional equations, one can cancel out the extraneous factors
and deduce that (see [Gr])
\[h_{K}^{-}=\frac{w_{K}}{2^{e}w_{K^{+}}}\lim_{s\rightarrow 0}
\frac{\zeta_{K}(s)}{\zeta_{K^{+}}(s)},\] where
$\frac{R_{K}}{R_{K^{+}}}=2^{e}$. But
\[\zeta_{K}(s)= \prod_{i=0}^{p-2}L(0,w^{i}),\;\; \zeta_{K^{+}}(s)=
\prod_{i \;\;even}^{p-2}L(0,w^{i}).\] Now observing that
$w_{K}=2p$ and $w_{K^{+}}=2$, we obtain the desired result.
\hspace*{1cm}$\square$

%----------
%----------------------

\subsection{Dirichlet $L$-values and Bernoulli numbers}

Recall that Bernoulli numbers $B_{n}$ are given by
\[\frac{t}{e^{t}-1}=\sum_{n=0}^{\infty}B_{n}\frac{t^{n}}{n!}.\]
Eg, $B_{0}=1, \; B_{1}=-\frac{1}{2},\; B_{2}= \frac{1}{6}$ etc.

\noin The $n$-th Bernoulli polynomial $B_{n}(X)$ is defined by
\[\frac{te^{Xt}}{e^{tX}-1}=\sum_{n=0}^{\infty}B_{n}(X)\frac{t^{n}}{n!}.\]
It is easy to see that \[B_{n}(X)=\sum_{i=o}^{n}{n \choose
i}B_{i}X^{n-i}.\] Eg, $B_{1}(X)=X-\frac{1}{2}$,
$B_{2}(X)=X^{2}-X+\frac{1}{6}$, etc.

\noin Now, for a Dirichlet character $\chi$ of conductor $f$, we
define the generalized Bernoulli numbers $B_{n,\chi}$ by
\[\sum_{a=1}^{f}\frac{\chi(a)te^{at}}{e^{ft}-1}=
\sum_{n=0}^{\infty}B_{n,\chi}\frac{t^{n}}{n!}.\] The following
well-known proposition allows us to express generalized Bernoulli
numbers in terms of Bernoulli polynomials (cf [Wa]).
\begin{prop}
If $g$ is any multiple of $f$, then
\[B_{n,\chi}=g^{n-1}\sum_{a=1}^{g}\chi(a)B_{n}\big(\frac{a}{g}\big).\]
\end{prop}
Proof: \begin{eqnarray*}
\sum_{n=0}^{\infty}g^{n-1}\sum_{a=1}^{g}\chi(a)B_{n}\big(\frac{a}{g}\big)\frac{t^{n}}{n!}
& = & \sum_{a=1}^{g}\chi(a)\frac{1}{g}\frac{(gt)e^{(\frac{a}{g})gt}}{e^{gt}-1}\\
& = &
\sum_{b=1}^{f}\sum_{c=0}^{h-1}\chi(b+cf)\frac{te^{(b+cf)t}}{e^{fht}-1}\;\;
\mx{ where }g=hf,\;\;a=b+cf\\
& = & \sum_{b=1}^{f}\frac{\chi(b)te^{bt}}{e^{ft}-1}\\
& = &\sum_{n=0}^{\infty}B_{n,\chi}\frac{t^{n}}{n!}. \;\;\;\;\qquad
\square
\end{eqnarray*}

\noin For example,
\begin{eqnarray*}
B_{1,\chi}& = & \sum_{a=1}^{f}\chi(a)(\frac{a}{f}-\frac{1}{2})
 =  \frac{1}{f}\sum_{a=1}^{f}\chi(a)(a-\frac{1}{2}f).\\
B_{2,\chi}& = & f\sum_{a=1}^{f}\chi(a)\Big(\frac{a}{f})^{2}-
\frac{1}{2}\frac{a}{f}+\frac{1}{6}\Big) =
\frac{1}{f}\sum_{a=1}^{f}\chi(a)\Big(a^{2}-fa+\frac{f^{2}}{6}\Big).
\end{eqnarray*}

\noin The generalized Bernoulli numbers can be relate to the
values of Dirichlet $L$-values as follows:
\begin{prop}
$L(1-n,\chi)=-\frac{B_{n,\chi}}{n}$, $n\geq 1$.
\end{prop}
For example, if $\chi$ is a Dirichlet character modulo $p$, we
have
\begin{eqnarray}\label{l1}L(0,\chi)= -B_{1,\chi}& = & -\frac{1}{p}
\sum_{n=1}^{p}\chi(n)\Big(n-\frac{1}{2}p\Big).\\
\label{l2}L(-1,\chi)=-B_{2,\chi}& = & -\frac{1}{2p}
\sum_{n=1}^{p}\chi(a)\Big(n^{2}-pn+\frac{p^{2}}{6}\Big).
\end{eqnarray}

%----------
%----------------------

\subsection{Some congruences involving Bernoulli numbers}

We require the following congruences involving Bernoulli numbers.

\begin{prop} (Kummer Congruence)
$\frac{B_{m}}{m}\equiv \frac{B_{n}}{n}$ if $m\equiv n \not\equiv
0$ mod $(p-1$).
\end{prop}

\noin Kummer's congruence can be proved in the following manner
(cf [B-S]): \\
let $g$ be a primitive root mod $p$. Consider
\begin{equation}
\label{bs1}
F(t)=\frac{gt}{e^{gt}-1}-\frac{t}{e^{t}-1}=\sum_{m=1}^{\infty}
(g^{m}-1)B_{m}\frac{t^{m}}{m!}.\end{equation} Letting $e^{t}-1=u$,
we can write
\[{F(t)=\frac{gt}{(1+u)^{g}-1}-\frac{t}{u}=tG(u),\mx { where }
G(u)=\frac{g}{(1+u)^{g}-1}-\frac{1}{u}=\sum_{k=1}^{\infty}c_{k}u^{k},\;\;c_{k}\in
\mbz}.\] Now,
\begin{equation}\label{bs2}G(u)=G(e^{t}-1)=\sum_{k=0}^{\infty}c_{k}(e^{t}-1)^{k}=
\sum_{m=1}^{\infty}{A_{m}}\frac{t^{m}}{m!}.\end{equation} But
$A_{m}$ are $p$-integral as they are integral linear combinations
of $c_{k}$'s. Further, they have period $(p-1)$ modulo $p$, as the
coefficients $r^{n}$ of $\frac{t^{n}}{n!}$ in $e^{rt}$ ($ r\geq
0$) have that periodicity by Fermat's little theorem
$r^{n+p-1}\equiv r^{n}$ modulo $p$. Comparing coefficients in
(\ref{bs1}) and (\ref{bs2}), we obtain
\[\frac{g^{m}-1}{m!}B_{m}=\frac{A_{m-1}}{(m-1)!}\;\Rightarrow\;
\frac{B_{m}}{m}(g^{m}-1)=A_{m-1}.\] If $p-1\not| m$, then
$g^{m}-1\not\equiv 0$ mod $p$ as $g$ is a primitive root mod $p$.
Clearly, $g^{m}-1$ has period $p-1$ mod $p$. Therefore,
$\frac{B_{m}}{m}$ also has period $p-1$ mod $p$ and is
$p$-integral.\hspace*{1cm}$\square$

\begin{prop}
$pB_{m}$ is $p$-integral, and $B_{m}$ is $p$-integral if
$(p-1)\not| m$.
\end{prop}
\begin{prop}
For an even integer $m$, $pB_{m}\equiv
\sum\limits_{a=1}^{p-1}a^{m}$ modulo $p^{2}$ if $p\geq 5$.
\end{prop}

\noin We can easily prove the above two propositions using the
following lemma.
\begin{lemma}
$(m+1)S_{m}(n)=\sum\limits_{k=0}^{m}{m+1 \choose
k}B_{k}n^{m+1-k}$, where $S_{m}(n)=1^{n}+2^{n}+\ldots+m^{n}$.
\end{lemma}
Proof:
\begin{eqnarray*}
\sum_{m=0}^{\infty}S_{m}(n)\frac{t^{m}}{m!} &=&\sum_{a=0}^{n-1}
\frac{e^{nt}-1}{e^{t}-1}=\frac{e^{nt}-1}{t}\frac{t}{e^{t}-1}
=\sum_{l=1}^{\infty}{n^{l}}\frac{t^{l-1}}{l!}
\sum_{k=0}^{\infty}B_{k}\frac{t^{k}}{k!}\\
\Rightarrow
\frac{S_{m}(n)}{m!}& = &\sum_{k=0}^{m+1}\frac{B_{k}}{(m+1-k)!k!}n^{m+1-k}\\
\Rightarrow (m+1)!\frac{S_{m}(n)}{m!}& = &\sum_{k=0}^{m+1}{m+1
\choose k}B_{k}n^{m+1-k}\hspace*{1cm} \square
\end{eqnarray*}

\noin In order to prove proposition 6.5, it is enough to show that
$pB_{m}\equiv S_{m}(p)$ modulo $p$. It is clear that
$S_{m}(p)\equiv 0$ mod $p$ if $(p-1)\not| m$ and $S_{m}(p)\equiv
p-1$ mod $p$ if $(p-1)| m$. By our lemma, we have
\begin{equation}\label{bs3}S_{m}(p)=pB_{m}+{m \choose 1}B_{m-1}\frac{p^{2}}{2}+{m
\choose 2}B_{m-2}\frac{p^{3}}{3}+\ldots+{m\choose
m}B_{0}\frac{p^{k+1}}{k+1}.\end{equation}

\noin Clearly, $\frac{p^{k+1}}{k+1}\equiv 0$ mod $p$ for $k\geq
2$, and $\frac{p^{k+1}}{k+1}$ is $p$-integral even for $k=1$.
Applying induction, let $pB_{j}$ be $p$-integral for $j<m$. Then,
$pB_{m}$ is $p$-integral as well, and we also obtain
$S_{m}(n)\equiv pB_{m}$ mod $p$ from (\ref{bs3}). Note that though
we need the result only for odd prime $p$, not that the above
proof works for $p=2$ as well, as $B_{n}$ vanishes for odd $n\geq
3$. $\square$\\

\noin To prove proposition 6.6, it suffices to establish that
$ord_{p}({m \choose k}B_{m-k}\frac{p^{k+1}}{k+1})\geq 2$ in view
of (\ref{bs3}). Since $pB_{m-k}$ is $p$-integral, we need only
$k-ord_{p}(k+1)\geq 2$. For $p\geq 5$ and $k\geq 2$, it is
obvious. For $k=1$, note that $B_{m-1}=0$ unless $m=2$, which
again follows trivially. $\square$

%----------
%----------------------

\subsection{A refined congruence for the Teichmuller character}

Let $w: (\mbz/p\mbz)^{\times}\kar \mu_{p-1}$ be the character
given by $w(n)\equiv n$ modulo $\wp$ where $\wp$ is any prime
ideal above $p$ in $\mbq(\mu_{p-1})$. The character $w$ is known
as the Teichmuller character. We have used the following
congruence for the Teichmuller character.
\begin{prop}
For $(n,p)=1$, we have $w(n)\equiv n^{p}$ modulo $\wp^{2}$ where
$\wp$ is a fixed prime above $p$ in $K=\mbq(\mu_{p-1})$.
\end{prop}
Proof: Let us recall Hensel's lemma:\\
 Let $R$ be a ring which is complete with respect to an ideal $I$
 and let $f(x)\in R[x]$. If
$f(a)\equiv 0$ mod ($f'(a)^{2}I$) then there exists $b\in R$ with
$b\equiv a$ modulo $(f'(a)I)$ such that $f(b)=0$. Further, $b$ is
unique if $f'(a)$ is a non-zero divisor in $R$.

Now let $K_{\wp}$ be the completion of $K$ at $\wp$. Let $R={\cal
O}_{\wp}$ be the completion of the ring of integers ${\cal O}$ of
$K$ with respect to $\wp$. Let $I=\wp^{2}$, then we can also think
of $R$ as the completion of ${\cal O}$ with respect to $I$.
Consider $f(x)=x^{p-1}-1$ and let $a=n^{p}$, where $(n,p)=1$.
Then, \[f(a)=(n^{p})^{p-1}-1\equiv 0\mx{ mod }\wp^{2}, \mx { as }
\#\Big(\frac{{\cal
O}_{\wp}}{\wp^{2}}\Big)^{\times}=\#\Big(\frac{{\cal
O}}{\wp^{2}}\Big)^{\times}=N\wp^{2}-N\wp=p(p-1).\] Moreover
$f'(a)=(p-1)a^{p-2}$ is not a zero-divisor in $R$. Therefore by
Hensel's lemma there exists a unique $b_{n}$ in $R$ such that
$b_{n}^{p-1}-1=0$ and $b_{n}\equiv n^{p}$ modulo $\wp^{2}$. Now,
if we define $w(n)=b_{n}$, we obtain the Teichmuller character
$w:\big(\frac{\mbz}{p\mbz}\big)^{\times}\kar \mu_{p-1}$ with the
more refined congruence $w(n)\equiv n^{p}$ modulo $\wp^{2}$.
\hspace*{1cm}$\square$

%----------
%----------------------

\subsection{Petersson inner product}

There is a measure on the upper half complex plane $\fh$ given by
$d\mu(\tau)=\frac{dx \;dy}{y^{2}}$ where $\tau=x+iy\in \fh$. It is
easy to show that $d\mu(\tau)$ is invariant under
$GL_{2}(\mbr)^{+}\subset Aut(\fh)$, i.e.,
$d\mu(\alpha\tau)=d\mu(\tau)$. In particular, the measure is
$SL_{2}(\mbz)$-invariant. As $\mbq\cup \{\infty\}$ is a countable
set of measure $0$, $d\mu$ suffices for integration over the
extended upper half plane $\fh^{*}=\fh\cup \mbq\cup \{\infty\}$.
Let $D^{*}$ be the fundamental domain for $SL_{2}(\mbz)$, i.e.,
\[D^{*}=\fh^{*}/SL_{2}(\mbz) =\{\tau \in \fh\mid
Re(\tau)\leq \frac{1}{2}, \;\;|\tau|\geq 1\}\cup \{\infty\}.\] For
a congruence subgroup $\Gamma$ of $SL_{2}(\mbz)$, we have $(\pm
I)\Gamma\ SL_{2}(Z)=\bigcup_{j}(\pm 1)\Gamma\alpha_{j}$ where $j$
runs over a finite set. Then, the fundamental domain for $\Gamma$
is given by
\[X(\Gamma)=\fh^{*}/\Gamma =\bigcup\alpha_{j}(D^{*}).\]
This allows us to integrate function of $\fh^{*}$ invariant under
$\Gamma$ by setting
\[\int\limits_{X(\Gamma)}\phi(\tau)d\mu(\tau)=
\int\limits_{\bigcup_{j}\alpha_{j}(D^{*})}\phi(\tau)d\mu(\tau)=
\sum_{j}\int\limits_{D^{*}}\phi(\alpha_{j}(\tau))d\mu(\tau).\] By
letting $V_{\Gamma}=\int\limits_{X(\Gamma)}d\mu(\tau)$, we can
define an inner product
\[<,>_{\Gamma}\;: \;S_{k}(\Gamma)\times M_{k}(\Gamma)\kar \mbc.\]
given by
\[<f,g>_{\Gamma}=\frac{1}{V_{\Gamma}}\int\limits_{X(\Gamma)}f(\tau)\overline{g(\tau)}
(Im(\tau))^{k}d\mu(\tau).\] Note that the integrand is invariant
under $\Gamma$. For the integral to converge, we need one of $f$
or $g$ to be a cusp form (see section 5.4 in [Di-S]). Clearly this
inner product is Hermitian and positive definite. When we take a
modular form $f\in M_{k}(\Gamma)- S_{k}(\Gamma)$, we can show that
$f$ is orthogonal
 under $<,>_{\Gamma}$ to all of $S_{k}(\Gamma)$. Thus,
we can think of the quotient space ${\cal
E}_{k}(\Gamma)=M_{k}(\Gamma)/S_{k}(\Gamma)$ as the complementary
subspace linearly disjoint from $S_{k}(\Gamma)$. This allows us to
write
\[S_{k}(\Gamma)=S_{k}(\Gamma)\oplus {\cal E}_{k}(\Gamma).\]

%----------
%----------------------

\subsection{Hecke operators}

For any $\alpha\in GL_{2}(\mbq)$, one can write the double coset
$\Gamma\alpha\Gamma=\bigcup_{i}\Gamma\alpha_{i}$ where
$\alpha_{i}$ runs over a finite set. We can define an action of
the double coset on $M_{k}(\Gamma)$ by setting $f|\Gamma \alpha
\Gamma=\sum f|[\alpha_{i}]$. It is easy to verify that these
operators preserve $M_{k}(\Gamma)$, $S_{k}(\Gamma)$ and ${\cal
E}_{k}(\Gamma)$.

We need to consider only the case $\Gamma=\Gamma_{1}(p)$. For any
integer $d$ such that $(d,p)=1$, we can define an operator $<d>$
as follows: we have a $ad-bp=1$ for some $a,\;b\in \mbz$. Taking
$\alpha= \left[\begin{array}{cc} a&b\\p&d\\\end{array}\right]\in
\Gamma_{0}(p)$, we obtain
\begin{eqnarray*}
<d> \;: \;M_{k}(\Gamma_{1}(p))\kar M_{k}(\Gamma_{1}(p)),\\
<d>f:=f|\Gamma_{1}(p)\alpha\Gamma_{1}(p)=f|[\alpha]_{k},\end{eqnarray*}
noting that $\Gamma_{1}(p)\alpha\Gamma_{1}(p)=\Gamma_{1}(p)\alpha$
as $\Gamma_{1}(p)$ is a normal subgroup of $\Gamma_{0}(p)$. The
operators $<d>$ are called diamond operators.

By taking $\alpha_{l}=\left[\begin{array}{cc}
1&0\\0&l\\\end{array}\right]$ for any prime $l$, we get an
operator $T_{l}=f|\Gamma \alpha_{l} \Gamma$ for any prime $l$. We
extend the definition of definition of Hecke operators to all
natural numbers inductively by setting \begin{eqnarray*}
T_{l^{r+1}}& = &T_{l}T_{l^{r}}- l^{k-1}<l> T_{l}^{r-1}\mx{ for
}r\geq 1.\\ T_{mn}& = & T_{m}T_{n}\;\;{\mx{when }} gcd(m,n)=1
\end{eqnarray*}

All these Hecke operators defined above are self adjoint with
respect to the Petersson inner product. For more details, see
chapter 5 of [Di-S]. A modular form is called an {\it eigenform}
if it is a simultaneous eigenform for all Hecke operators $T_{n}$
and $<d>$, $(d,p)=1$.

%---------------
%---------------------

\subsection{Ingredients from commutative algebra}

The results proved below are required for the lifting lemma of
Deligne and Serre in section 4.1.

\subsubsection{Minimal primes}
Let $A$ be a commutative ring with 1. A prime ideal $\wp$ of $A$
is called a {\it minimal prime} if it the smallest prime ideal
(containing $0$) in $A$. Such a prime exists by Zorn's lemma on
the (non-empty as $1\in A$) set $S$ of primes ideals of $A$ with
the partial order $I\leq J$ when $J\subset I$, noting that any
descending chain in $S$ has its intersection as an upper bound in
$S$.
\begin{prop}
A minimal prime $\wp$ of $A$ is contained in the set $Z$ of
zero-divisors of $A$.
\end{prop}
Proof: Note that $x,\;y\in D=A-Z \Rightarrow xy\in D$. Thus $D$ is
a multiplicative set. On the other hand, $S=A-\wp$ is a maximal
multiplicative closed set (as $\wp$ is a minimal prime). If
$D\not\subset S$, then $SD$ would be a multiplicative set strictly
larger than $S$. Therefore, $D\subset S$ and $\wp\subset Z$.
\hspace*{1cm}$\square$

\subsubsection{Associated primes and support primes}
Let $A$ be a
commutative ring and $M$ be an $A$-module. The annihilator of a
submodule $N$ of $M$ is defined as
\[Ann_{A}(N)=\{a\in A|an=0 \;\forall n\in N\}.\] Clearly,
$Ann_{A}(N)$ is an ideal of $A$. For an element $m\in M$, we can
define its annihilator as $Ann_{A}(m)= \{a\in A|am=0\}$.
\begin{defn} A prime ideal $\wp$ of $A$ is called an
{\bf associated prime} if $\wp$ is the annihilator of some element
of $M$. The set of associated primes of $M$ in $A$ is denoted by
$Assoc_{A}(M)$.
\end{defn}
\begin{prop} If $M$ is non-zero and $A$ is Noetherian, then
$Assoc_{A}(M)$ is non-empty.
\end{prop}
Proof: Consider the set $S$ of ideals ($\neq A$) of $A$ which are
annihilators of some element of $M$. As $A$ is Noetherian, $S$ has
a maximal element, say $\wp$, which is necessarily the annihilator
of some element $m$ in $M$. Let $x,\;y\in A$ such that $xy\in \wp$
but $y\not\in \wp$. Then $ym\neq 0$, but $\wp \subset
(\wp,x)\subset Ann_{A}(ym)\in S$. It follows that
$Ann_{A}(ym)=(\wp,x)=\wp$ by maximality of $\wp$. Therefore $x\in
\wp$, and hence $\wp$ is an associated prime. \hspace*{1cm}
$\square$
\begin{defn}
 A prime ideal $\wp$ of $A$ is called a {\bf support
 prime}
of $M$ if $M_\wp\neq 0$.
\end{defn}
The set of support primes of $M$ in $A$ is denoted by
$Supp_{A}(M)$.
\begin{prop}
Let $A$ be Noetherian and $M$ be a finitely generated $A$-module.
Then $\wp\in Supp_{A}(M)\Leftrightarrow Ann_{A}(M)\subset\wp$
\end{prop}
Proof: Let $Ann_{A}(M)\not\subset \wp$. Then there exists $s\in
A-\wp$ such that $sM=0$, hence $M_{\wp}=0$. Contra-positively,
$\wp\in Supp_{A}(M)$ implies $Ann_{A}(M)\subset \wp$.

For the converse, let $m_{1},\ldots , m_{r}$ generate $M$ as an
$A$-module. If $M_{\wp}=0$, then we can find $s_{i}\in A-\wp$ such
that $s_{i}m_{i}=0$. Now $s=s_{1}\ldots s_{r}\in A-\wp$
annihilates $M$, hence $Ann_{A}(M)\not\subset \wp$.
\hspace*{1cm}$\square$

\begin{prop}
$Assoc_{A}(M)\subset Supp_{A}(M)$.
\end{prop}
Proof: Let $\wp$ be an associated prime of $M$, say
$\wp=Ann_{A}(m)$ for some $m\in M$. If $M_{\wp}=0$ then there
exists $s\in A- \wp$ such that $sm=0$. But it would mean $s\in
Ann_{A}(m)=\wp$, which is a contradiction. Thus, $M_{\wp}\neq 0$
and $\wp$ must be a support prime of $M$. \hspace*{1cm}$\square$
\begin{prop}
Let $A$ be a Noetherian ring and $\wp$ be a support prime. Then
$\wp$ contains an associated prime $\fq$ of $M$.
\end{prop}
Proof: If $\wp$ is a support prime, $M_{\wp}\neq 0$. Then there
must exist some $x\in M$ such that $(Ax)_{\wp}\neq 0$. Thus, there
exists an associated prime $\fq$ of the $A$-module $(Ax)_{\wp}$.
Hence there is an element $0\neq \frac{y}{s}$ of $(Ax)_{\wp}$ with
$y\in Ax$ and $s\not\in \wp$ such that $\fq$ is the annihilator of
$\frac{y}{s}$. Now, if there exists $b\in \fq-\wp$, then
$b\frac{y}{s}=0$ would imply $\frac{y}{s}=0$, which is a
contradiction.

Now we still have to show that $\fq$ is an associated prime of $M$
as well. Let $b_{1},\;\ldots b_{n}$ be a set of generators of
$\fq$. Then, there exists $t_{i}\in A-\wp$ such that
$b_{i}t_{i}y=0$. Let $t=t_{1}.\ldots . t_{n}$. Then, $\fq$ is the
annihilator of $ty\in M$. \hspace*{1cm}$\square$
\begin{cor}
If $\wp$ is a minimal prime in the support of $M$, then $\wp$ is
also an associated prime when $A$ is Noetherian.
\end{cor}
Proof: As $\wp$ must contain an associated prime, we get our
result by minimality of $\wp$. $\square$


\begin{thebibliography}{C-W 2}

\bibitem[BS]{BS}
Borevich, Z. I., Shafarevich, I. R.; {\itshape Number Theory},
Academic Press, 1966.

\bibitem[C]{C}
Carlitz, L.; {\it A generalization of Maillet's determinant and a
bound for the first factor of the class number}, Proc. A.M.S. 12,
256--261, 1961.

\bibitem[C-O]{C-O}
Carlitz, L. Olson F.R.; {\it Maillet's determinant}, Proc. A.M.S.
6, 265--269, 1955.

\bibitem[D]{D}
Dalawat, C.S.; {\it Ribet's modular construction of unramified
$p$-extensions of $\mbq(\mu_{p})$} (to appear)


\bibitem[D-S]{D-S}
Deligne P., Serre, J-P.; {\it Formes modulaires de poids 1}, Ann.
Scient. Ec. Norm. Sup., $4^{e}$ serie, 7, 507--530, 1974.

\bibitem[Di-S]{Di-S}
Diamond, F., Shurman S.; {\itshape A First Course on Modular
Forms}, Springer, 2005.

\bibitem[Gr]{Gr}
Greenberg, R.; {\itshape A generalization of Kummer's criterion},
Inventiones Math. 21, 247--254, 1973.
\bibitem[La]{La}
Lang, S.; {\itshape Algebra}, Springer-Verlag, 2002.

\bibitem[R]{R}
Ribet, K.; {\itshape A modular construction of unramified
$p$-extensions of $\mbq(\mu_{p})$}, Inventiones Math. 34,
151--162, 1976.

\bibitem[Wa]{Wa}
Washington, L.; {\itshape Introduction to Cyclotomic Fields},
Springer-Verlag, 1997.
\end{thebibliography}
\end{document}